\documentclass[12pt]{amsart}
\usepackage{amsmath, amsfonts, amssymb, amsthm,hyperref,mathtools,array}
\hypersetup{hypertex=true,
	colorlinks=true,
	linkcolor=blue,
	anchorcolor=blue,
	citecolor=blue}
\usepackage{bm}
\allowdisplaybreaks[4]
\def\ord{{\rm ord}}

\def\x{{\bm x}}

\def\0{{\bm 0}}

\def\pmod #1{\ ({\rm{mod}}\ #1)}
\def\mod #1{\ {\rm mod}\ #1}
\def\bN{{\bf N}}
\def\bH{{\bf H}}
\def\GL{{\rm GL}}

\theoremstyle{plain}
\newtheorem{theorem}{Theorem}[section]
\newtheorem{lemma}{Lemma}

\newtheorem{conjecture}{Conjecture}

\newtheorem{definition}{Definition}

\theoremstyle{definition}

\theoremstyle{remark}

\makeatletter
\@namedef{subjclassname@2020}{%
	\textup{2020} Mathematics Subject Classification}
\makeatother
\vspace{4mm}

\begin{document}
	
	\title[On the Burgess bound and Pythagorean triples involving primitive roots] 
	{On the Burgess bound and Pythagorean triples involving primitive roots}
	\author[H.-L. Wu and H.-X. Ni]{Hai-Liang Wu and He-Xia Ni*}
	
	\address {(Hai-Liang Wu) School of Science, Nanjing University of Posts and Telecommunications, Nanjing 210023, People's Republic of China}
	\email{\tt whl.math@smail.nju.edu.cn}
	
	\address {(He-Xia Ni) Department of Applied Mathematics, Nanjing Audit University, Nanjing 211815, People's Republic of China}
	\email{\tt nihexia@yeah.net}

	\keywords{Pythagorean triples, primitive roots, finite fields.
		\newline \indent This research was supported by the Natural Science Foundation of China (Grant Nos. 12101321 and 12371004) and the Natural Science Foundation of the Higher Education Institutions of Jiangsu Province (Grant No. 25KJB110010).
		\newline \indent *Corresponding author.}
	
	\begin{abstract}
	By combining the recent results of Pierce and Xu on the higher-dimensional Burgess bound with the classical one-dimensional Burgess bound, we prove that for any sufficiently large prime $p$, there exists a Pythagorean triple $(a_p,b_p,c_p)$ with $a_p,b_p,c_p\in\mathbb{Z}\cap(0,p)$ such that $a_pb_p/2$ is a primitive root modulo $p$. This confirms a conjecture of Z.-W. Sun for all sufficiently large primes. 
	\end{abstract}
	\maketitle
	
	\section{Introduction}
	\setcounter{lemma}{0}
	\setcounter{theorem}{0}
	\setcounter{equation}{0}
	\setcounter{conjecture}{0}
	\setcounter{remark}{0}
	\setcounter{corollary}{0}

	\subsection{Notation} Throughout this paper, $p$ denotes a prime and let $\mathbb{F}_p=\mathbb{Z}/p\mathbb{Z}$ be the finite field with $p$ elements. Let $\mathbb{F}_p^*=\mathbb{F}_p\setminus\{0\}$ be the multiplicative cyclic group of all non-zero elements over $\mathbb{F}_p$.  An element $g\in\mathbb{F}_p^*$ is said to be a primitive element of $\mathbb{F}_p$ if $g$ is a generator of the cyclic group $\mathbb{F}_p^*$. The set of all primitive elements of $\mathbb{F}_p$ is written as $\mathcal{P}_p$. In addition, given any integer $m$, we say that $m$ is a primitive root modulo $p$ if $m\mod{p\mathbb{Z}}\in\mathbb{F}_p$ is primitive. Let $\widehat{\mathbb{F}_p^*}$ be the group of all multiplicative characters of $\mathbb{F}_p$. For any character $\psi: \mathbb{F}_p^*\rightarrow\mathbb{C}^*$, we additionally define $\psi(0)=0$ and $\ord(\psi)$ is the order of $\psi$. In addition, we use $\chi_0$ to denote the trivial character, that is, $\chi_0(x)=1$ for any $x\in\mathbb{F}_p^*$ and $\chi_0(0)=0$. For simplicity, the Legendre symbol $(\frac{\cdot}{p})$ is abbreviated as $\rho$, i.e., $\rho$ is a character with $\ord(\rho)=2$ defined by 
   $$\rho(x)=\begin{cases}
   	1  & \mbox{if}\ x\in\{y^2: y\in\mathbb{F}_p^*\},\\
   	0 & \mbox{if}\ x=0,\\
   	-1 & \mbox{otherwise.}
   \end{cases}$$
	
	Also, we use the symbol $\#S$ to denote the cardinality of a set $S$ and let  $$\GL_n(\mathbb{F}_p)=\left\{M: M\ \text{is an $n\times n$ invertible matrix over $\mathbb{F}_p$}\right\}.$$
	
	Let $f$ and $g$ be functions defined on the set of all primes. As usual, $f(p)\ll g(p)$ means that there exist a positive real number $C$ and a real number $x_0$ such that 
	$$|f(p)|\le C\cdot g(p)$$
	for any prime $p\ge x_0$. The symbol $f(p)\asymp g(p)$ means that $f(p)\ll g(p)$ and $g(p)\ll f(p)$. If $g(p)$ is non-vanishing for any prime $p\ge x_0$, then the symbol $f(p)=o(g(p))$ means that 
	$$\lim_{p\rightarrow+\infty}\frac{f(p)}{g(p)}=0.$$

	\subsection{Background and motivation}
	
	Let $a,b,c\in\mathbb{Z}^+$. We say that $(a,b,c)$ is a Pythagorean triple if $a^2+b^2=c^2$. Geometrically, $(a,b,c)$ forms a right triangle with $c$ as its hypotenuse and $a, b$ as its legs. Pythagorean triples have been studied since ancient times, with notable contributions by Pythagoras, Euclid, Diophantus, and many other mathematicians. Pythagorean triples are also closely related to number theory. For example, it is well known that each positive integer solution to the Pythagorean equation
	\begin{equation}\label{Eq. Pythagorean equation}
		x^2+y^2=z^2
	\end{equation}
	is of the form 
	$$\left(d(m^2-n^2), d(2mn), d(m^2+n^2)\right)$$
	where $d,m,n\in\mathbb{Z}^+$ with $m>n$, $\gcd(m,n)=1$ and $m\not\equiv n\pmod 2$. On the other hand, the classical congruent number problem asks: given a positive integer $n$, can one find a right triangle with all three sides rational such that its area is precisely $n$? This is a very interesting and profound problem, which has deep connections with elliptic curves and BSD conjecture.
	
	Inspired by the above research history, Z.-W. Sun initiated the study of the problem: given a prime $p$, does there exist a positive integer solution $(x_p,y_p,z_p)$ to (\ref{Eq. Pythagorean equation}) that satisfies certain algebraic property over $\mathbb{F}_p$? For example, Z.-W. Sun \cite{Sun15} posed the following conjecture involving quadratic residues modulo $p$.
	
	\begin{conjecture}[Z.-W. Sun]\label{Conj. Sun quadratic}
		Let $p>100$ be a prime. Then, for any $\varepsilon_1,\varepsilon_2,\varepsilon_3\in\{-1,1\}$, there is a positive integer solution $(x_p,y_p,z_p)$ to (\ref{Eq. Pythagorean equation}) such that $0<x_p,y_p,z_p<p$ and 
		$$\left(\rho(x_p),\rho(y_p),\rho(z_p)\right)=\left(\varepsilon_1,\varepsilon_2,\varepsilon_3\right).$$
	\end{conjecture}
	By applying the higher-dimensional Burgess bound obtained by Pierce and Xu \cite{PX}, Xi and Zheng \cite{XZ} proved that this conjecture holds for all sufficiently large primes. 
	
	On the other hand, primitive elements also play an important role in the theory of finite fields and possess a variety of interesting combinatorial properties. For example, given an integer $n\ge3$, Cohen and his collaborators \cite{Cohen15} showed that if the odd prime power $q$ is large enough and the characteristic of the finite field $\mathbb{F}_q$ is greater than $n$, then there is an element $x$ such that the $n$ consecutive elements 
	$$x,x+1,\cdots,x+n-1$$
	are all primitive elements of $\mathbb{F}_q$. For more results on this topic, readers may refer to \cite{Cohen03,Cohen10,Cohen15,Cohen21}. Along this line, Z.-W. Sun \cite{Sun19} posed the following conjecture related to Pythagorean triples and primitive roots.
	
	\begin{conjecture}[Z.-W. Sun] \label{Conj. primitive}
		Let $p>7$ be a prime. Then, there exists a right triangle with all side lengths in $\{1,2,\dots,p-1\}$ such that its area is precisely a primitive root modulo $p$.
	\end{conjecture}
	
	\subsection{Main results} By combining the recent results of Pierce and Xu on the higher-dimensional Burgess bound with the classical one-dimensional Burgess bound, we can confirm this conjecture for all sufficiently large primes. Now we state our main theorem.
	
	\begin{theorem}\label{Thm. A}
		Let $p$ be a sufficiently large prime. Then, there exists a Pythagorean triple $(a_p,b_p,c_p)$ with $a_p,b_p,c_p\in\mathbb{Z}\cap(0,p)$ such that $a_pb_p/2$ is a primitive root modulo $p$.
	\end{theorem}
	
	\subsection{Outline of the paper} In Section 2, we will introduce some necessary results on the Burgess bound. The proof of our theorem shall be given in Section 3.

	\section{The Burgess bound}
	\setcounter{lemma}{0}
	\setcounter{theorem}{0}
	\setcounter{equation}{0}
	\setcounter{conjecture}{0}
	\setcounter{remark}{0}
	\setcounter{corollary}{0}
	
	In 1957, Burgess \cite{Burgess57} obtained the following celebrated result on character sums over a short interval, which will be used in the proof of Theorem \ref{Thm. A}.
	
	\begin{theorem}[Burgess]\label{Thm. Burgess}
		Let $p$ be a prime, $N\in\mathbb{R}$ and $H\in\mathbb{R}_{\ge 1}$. For any non-trivial character $\chi\in\widehat{\mathbb{F}_p^*}$ and positive integer $r$, we have 
		$$\left|\sum_{N<x\le N+H}\chi(x)\right|\ll_r H^{1-\frac{1}{r}}\cdot p^{\frac{r+1}{4r^2}}\cdot \log p.$$
	\end{theorem}
	
	In 2020, Pierce and Xu \cite{PX} investigated the higher-dimensional Burgess bound. Let the integer $n\ge2$ and $F(x_1,\cdots,x_n)\in\mathbb{Z}[x_1,\cdots,x_n]$ be a homogeneous polynomial of degree $D$. Let 
	$$\bN=(N_1,\cdots, N_n)\in\mathbb{R}^n$$ 
	and 
	$$\bH=(H_1,\cdots,H_n)\in\mathbb{R}_{\ge1}^n.$$
	For any non-trivial $\chi\in\widehat{\mathbb{F}_p^*}$, define the character sum
	$$S(F; \bN,\bH)=\sum_{\substack{\x=(x_1,\cdots,x_n)\in\mathbb{Z}^n\\ x_i\in(N_i,N_i+H_i]}}\chi(F(\x)).$$
	
	To state the result of Pierce and Xu, we need to introduce the following definition \cite[Condition 1.1]{PX}.
	
	\begin{definition}\label{Def. admissible polynomial}
     Let $p$ be a prime and $\Delta\in\mathbb{Z}^+$. Let $f\in\mathbb{F}_p[x_1,\cdots,x_n]$ with the decomposition $f=g^{\Delta}h$, where $g,h\in\mathbb{F}_p[x_1,\cdots,x_n]$ and $h$ is $\Delta$-th power-free over $\mathbb{F}_p$. If $h$ has the property that it cannot be made independent of (at least) one variable after an invertible linear transformation, i.e. there exists no linear change of variables $A\in \GL_n(\mathbb{F}_p)$ such that $h(\x A)\in\mathbb{F}_p[x_2,\cdots,x_n]$, then we say that $f$ is $(\Delta,p)$-admissible. 
	\end{definition}
	
	For any $\bH=(H_1,\cdots,H_n)\in\mathbb{R}_{\ge1}^n$, let $\|\bH\|=H_1\cdots H_n$, 
	$$H_{\max}=\max\{H_i: 1\le i\le n\}\ \text{and}\ H_{\min}=\min\{H_i: 1\le i\le n\}.$$
	Now we introduce the result of Pierce and Xu \cite[Theorem 1.1]{PX}, which will play an important role in the proof of Theorem \ref{Thm. A}. 
	
	\begin{theorem}[Pierce and Xu]\label{Thm. PX on Burgess bound}
		Let $p$ be a prime and $\chi\in\widehat{\mathbb{F}_p^*}$ with $\ord(\chi)=\Delta>1$. Fix an integer $n\ge2$. For any $r\in\mathbb{Z}^+$, define
		$$\Theta_{n,r}=\left\lfloor\frac{r-1}{n-1}\right\rfloor.$$
		Let $\bH=(H_1,\cdots,H_n)\in \mathbb{R}_{\ge1}^n$ with $H_{\max}H_{\min}<p^{1+1/(2\Theta_{n,r})}$. Let  $F(x_1,\cdots,x_n)\in\mathbb{Z}[x_1,\cdots,x_n]$ be a homogeneous polynomial of degree $D$ such that its reduction modulo $p$ is $(\Delta,p)$-admissible. Then, for any $\bN=(N_1,\cdots, N_n)\in\mathbb{R}^n$ and any $r\in\mathbb{Z}^+$ we have 
		$$\left|S(F; \bN,\bH)\right|\ll \|\bH\|^{1-\frac{1}{2r}}\cdot H_{\min}^{-\frac{1}{2r}}\cdot p^{\frac{1+n\Theta_{n,r}}{4r\Theta_{n,r}}}\cdot \left(\log p\right)^{n+1},$$
		in which the implied constant depends only on $D,\Delta, n,r$ and is otherwise independent of $F$. 
	\end{theorem}

	\section{Proof of Theorem \ref{Thm. A}}
	\setcounter{lemma}{0}
	\setcounter{theorem}{0}
	\setcounter{equation}{0}
	\setcounter{conjecture}{0}
	\setcounter{remark}{0}
	\setcounter{corollary}{0}
	
	We begin with the known result on the characteristic function $1_{\mathcal{P}_p}$ of $\mathcal{P}_p$ (cf. \cite{Cohen03,Cohen10,Cohen15,Cohen21}). 
	
	\begin{lemma}\label{Lem. characteristic function of P}
		Let $p$ be a prime and $\mathcal{P}_p$ be the set of all primitive elements of $\mathbb{F}_p$. Then 
		$$1_{\mathcal{P}_p}(x)=\theta_{p-1}\sum_{d\mid p-1}\frac{\mu(d)}{\varphi(d)}\sum_{\substack{\chi\in\widehat{\mathbb{F}_p^*} \\ \ord(\chi)=d}} \chi(x)=\theta_{p-1}\sum_{\chi\in\widehat{\mathbb{F}_p^*}}c_{\chi}\cdot \chi(x)
		=\begin{cases}
			1  & \mbox{if}\ x\in\mathcal{P}_p,\\
			0 & \mbox{otherwise},
		\end{cases}$$
		where $\mu$ is the M\"obius function, $\varphi$ is the Euler totient function, $\theta_{p-1}=\frac{\varphi(p-1)}{(p-1)}$, and $c_{\chi}=\frac{\mu(\ord(\chi))}{\varphi(\ord(\chi))}$ for any $\chi\in\widehat{\mathbb{F}_p^*}$. 
	\end{lemma}
	
		For any positive integer $n$, let 
	$$\omega_n=\#\left\{p: p\mid n\ \text{and}\ p\ \text{is a prime}\right\}$$
	be the number of all distinct prime divisors of $n$. We need the following result due to Robin \cite[Theorem 11]{Robin}.
	
	\begin{lemma}\label{Lem. bound due to Robin}
		For any integer $n\ge 3$, we have 
		$$\omega_n<1.3841\frac{\log n}{\log\log n}.$$
	\end{lemma}
	
		Let 
		$$W_{p-1}=\#\left\{d\in\mathbb{Z}^+: d\mid p-1\ \text{and}\ d\ \text{is square-free}\right\}$$
		be the number of positive square-free divisors of $p-1$. Applying Lemma \ref{Lem. bound due to Robin}, we obtain 
		$$W_{p-1}=2^{\omega_{p-1}}=e^{\omega_{p-1}\cdot \log 2}\le (p-1)^{\frac{1.3841\cdot\log 2}{\log\log(p-1)}}$$
		for any prime $p\ge 5$. From this, we immediately obtain the following result.
		
		\begin{lemma}\label{Lem bound for W}
			For any positive real number $\varepsilon$, if $p$ is large enough, then  
				$$W_{p-1}\ll_{\varepsilon} p^{\varepsilon}.$$
		\end{lemma}

	Now we are in a position to prove our theorem.
	
	{\noindent\bfseries Proof of Theorem \ref{Thm. A}}. Suppose that $p$ is a sufficiently large prime. Let $H=\lfloor p^{17/48} \rfloor$ and 
	$$X_p=\left\{(x,y)\in\mathbb{Z}\times\mathbb{Z}: 2H<x\le 3H, H<y\le 2H\right\}.$$
	Define 
	$$a(x,y)=x^2-y^2,\  b(x,y)=2xy,\  c(x,y)=x^2+y^2.$$
	Since $p$ is large enough, it follows from the definition of $X_p$ that 
	\begin{equation}\label{Eq. a,b,c are positive integers less than p}
		0<a(x,y),b(x,y),c(x,y)\le 13H^2<p
	\end{equation}
	for any $(x,y)\in X_p$. 
	
	Recall that $\rho\in\widehat{\mathbb{F}_p^*}$ with $\ord(\rho)=2$. Next we consider the sum 
	$$\sum_{(x,y)\in X_p}\rho\left(a(x,y)b(x,y)/2\right)=\sum_{(x,y)\in X_p}\rho\left((x^2-y^2)xy\right).$$
	Let $F(\x)=F(x_1,x_2)=(x_1^2-x_2^2)x_1x_2$ and $f(\x)$ be its reduction modulo $p$.  As 
	$$x_1-x_2,\ x_1+x_2,\ x_1,\ x_2$$
	are mutually distinct polynomials over $\mathbb{F}_p$, the polynomial $f(\x)$ is square-free over $\mathbb{F}_p$. Suppose that there exists a matrix $M\in\GL_2(\mathbb{F}_p)$ such that $f(\x M)$ is independent of $x_2$. As $f(\x)$ is a homogeneous polynomial of degree $4$, we have 
	$$f(\x M)=\lambda x_1^4$$
	for some $\lambda\in\mathbb{F}_p^*$. This implies that 
	$$f(x_1,x_2)=\lambda g(x_1,x_2)^4$$
	for some $g(x_1,x_2)\in\mathbb{F}_p[x_1,x_2]$ with degree $1$, which contradicts the fact that $f(\x)$ is square-free over $\mathbb{F}_p$. Thus, by Definition \ref*{Def. admissible polynomial} the polynomial $f(\x)$ is $(2,p)$-admissible. Let $\bN=(2H, H)$ and $\bH=(H,H)$. Then $\bH$ clearly satisfies the required condition of Theorem \ref{Thm. PX on Burgess bound}. Applying Theorem \ref{Thm. PX on Burgess bound} with $n=2$ and $r=10$, one can verify that 
	\begin{equation}\label{Eq. sum involving quadratic character}
		\left|\sum_{(x,y)\in X_p}\rho\left((x^2-y^2)xy\right)\right|=\left|S(F; \bN, \bH)\right|\ll H^{\frac{37}{20}}p^{\frac{19}{360}}(\log p)^3.
	\end{equation}
	Noting that $H^2\asymp p^{17/24}$, we obtain $H^{\frac{-3}{20}}p^{\frac{19}{360}}\asymp p^{\frac{-1}{2880}}$. Applying this to (\ref{Eq. sum involving quadratic character}), we have 
	\begin{align}\label{Eq. sum of quadratic character is o(H2)}
		\left|\sum_{(x,y)\in X_p}\rho\left((x^2-y^2)xy\right)\right|
&\ll H^2H^{\frac{-3}{20}}p^{\frac{19}{360}}(\log p)^3\notag\\
&\ll H^2p^{\frac{-1}{2880}}(\log p)^3\notag\\
&=o(H^2).
	\end{align}
	Assembling (\ref{Eq. a,b,c are positive integers less than p}) and (\ref{Eq. sum of quadratic character is o(H2)}) gives 
	\begin{align*}
		2\cdot \#\left\{(x,y)\in X_p: \rho\left((x^2-y^2)xy\right)=-1\right\}
&=\sum_{(x,y)\in X_p}\left(1-\rho\left((x^2-y^2)xy\right)\right)\\
&=\# X_p-\sum_{(x,y)\in X_p}\rho\left((x^2-y^2)xy\right)\\
&=H^2+o(H^2)\\
&>0.
	\end{align*}
	Hence, there is an element $(x_0,y_0)\in X_p$ such that $\rho((x_0^2-y_0^2)x_0y_0)=-1$. Let 
	\begin{equation}\label{Eq. a0, b0, c0}
		a_0=a(x_0,y_0)=x_0^2-y_0^2,\ b_0=b(x_0,y_0)=2x_0y_0,\ c_0=c(x_0,y_0)=x_0^2+y_0^2,
	\end{equation}
	and $L=\lfloor(p-1)/c_0\rfloor$. Noting that $0<c_0\le 13H^2$ and $H^2\asymp p^{17/24}$, we obtain 
	\begin{equation}\label{Eq. L is large}
		L\gg \frac{p}{H^2}\gg p^{\frac{7}{24}}.
	\end{equation} 
	
	Next we focus on 
	$$N_p=\#\left\{1\le k\le L: (ka_0)(kb_0)/2=k^2(x_0^2-y_0^2)x_0y_0\ \text{is a primitive root modulo $p$}\right\}.$$
	Letting $A_0=(x_0^2-y_0^2)x_0y_0$, by Lemma \ref{Lem. characteristic function of P} one can verify that 
	\begin{align}\label{Eq. formula of Np}
		N_p
&=\sum_{1\le k\le L} 1_{\mathcal{P}_p}(k^2 A_0)\notag\\
&=\theta_{p-1}\sum_{1\le k\le L}\sum_{\chi\in\widehat{\mathbb{F}_p^*}}c_{\chi}\cdot \chi(A_0)\cdot \chi(k^2)\notag\\
&=\theta_{p-1}\sum_{\chi\in\widehat{\mathbb{F}_p^*}}c_{\chi}\cdot \chi(A_0)\cdot S(\chi),
	\end{align}
	where 
	$$S(\chi)=\sum_{1\le k\le L}\chi(k^2)=\sum_{1\le k\le L}\chi^2(k).$$
	Suppose that $\chi$ is neither $\chi_0$ nor $\rho$. Since $\chi^2\neq \chi_0$, by applying Theorem \ref{Thm. Burgess} with $r=7$, one can verify that 
	\begin{equation}\label{Eq. S(chi) when chi is neither trivial nor quadratic}
	\left|S(\chi)\right|\ll L^{\frac{6}{7}}\cdot p^{\frac{2}{49}}\cdot \log p.
	\end{equation}
	On the other hand, observe that  
	$$\#\left\{\chi\in\widehat{\mathbb{F}_p^*}: \ord(\chi)=d\right\}=\varphi(d)$$
	for any positive divisor $d$ of $p-1$. Hence, by Lemma \ref{Lem bound for W} we have 
	\begin{align}\label{Eq. sum of coeffient c chi}
		\sum_{\chi\in\widehat{\mathbb{F}_p^*}}\left|c_{\chi}\right|
&=\sum_{\chi\in\widehat{\mathbb{F}_p^*}}\frac{|\mu(\ord(\chi))|}{\varphi(\ord(\chi))}\notag\\
&=\sum_{d\mid p-1}\frac{|\mu(d)|}{\varphi(d)}\#\left\{\chi\in\widehat{\mathbb{F}_p^*}: \ord(\chi)=d\right\}\notag\\
&=\sum_{d\mid p-1}|\mu(d)|\notag\\
&=W_{p-1}\notag\\
&\ll_{\varepsilon} p^{\varepsilon},
	\end{align}
where $\varepsilon$ is an arbitrary positive real number.  Combining (\ref{Eq. L is large}) and (\ref{Eq. S(chi) when chi is neither trivial nor quadratic}) with (\ref{Eq. sum of coeffient c chi}) and noting that $|\chi(A_0)|=1$, for any $\varepsilon\in(0, 1/1176)$ one can verify that 
	\begin{align}\label{Eq. the error term}
		\left|\sum_{\chi\in\widehat{\mathbb{F}_p^*}\setminus\{\chi_0,\rho\}}c_{\chi}\cdot \chi(A_0)\cdot S(\chi)\right|
&\le L^{\frac{6}{7}}\cdot p^{\frac{2}{49}}\cdot \log p\sum_{\chi\in\widehat{\mathbb{F}_p^*}}\left|c_{\chi}\right|\notag\\
&\ll_{\varepsilon} LL^{\frac{-1}{7}}\cdot p^{\frac{2}{49}}\cdot \log p\cdot p^{\varepsilon}\notag\\
&\ll_{\varepsilon} L\cdot p^{\frac{-1}{24}}\cdot p^{\frac{2}{49}}\cdot \log p\cdot p^{\varepsilon}\notag\\
&\ll_{\varepsilon} L\cdot p^{\frac{-1}{1176}+\varepsilon}\cdot \log p\notag\\
&=o(L).
	\end{align}
	
When $\chi$ is either $\chi_0$ or $\rho$, it is clear that 
	\begin{equation}\label{Eq. S(chi) when chi is trivial or quadratic}
		S(\chi)=L.
	\end{equation}
	
	Applying (\ref{Eq. the error term}) and (\ref{Eq. S(chi) when chi is trivial or quadratic}) to (\ref{Eq. formula of Np}) and noting that $\rho(A_0)=-1$, we obtain 
	\begin{align*}
		\frac{1}{\theta_{p-1}}\cdot N_p
&=c_{\chi_0}\cdot \chi_0(A_0)\cdot S(\chi_0)+c_{\rho}\cdot \rho(A_0)\cdot S(\rho)+\sum_{\chi\in\widehat{\mathbb{F}_p^*}\setminus\{\chi_0,\rho\}}c_{\chi}\cdot \chi(A_0)\cdot S(\chi)\\
&\ge 2L-\left|\sum_{\chi\in\widehat{\mathbb{F}_p^*}\setminus\{\chi_0,\rho\}}c_{\chi}\cdot \chi(A_0)\cdot S(\chi)\right|\\
&=2L-o(L)\\
&>0.
	\end{align*}
	Hence, there exists an integer $k\in[1,L]$ such that $k^2A_0=k^2(x_0^2-y_0^2)x_0y_0$ is a primitive root modulo $p$.  Noting that $(a_0,b_0,c_0)$ defined by (\ref{Eq. a0, b0, c0}) is a Pythagorean triple with $a_0,b_0,c_0\in\mathbb{Z}\cap(0, p)$ and $L=\lfloor(p-1)/c_0\rfloor$, we immediately see that 
   $$(ka_0,kb_0,kc_0)$$
	is also a Pythagorean triple with $ka_0,kb_0,kc_0\in\mathbb{Z}\cap(0, p)$ and that $(ka_0)(kb_0)/2$ is a primitive root modulo $p$.
	
	In view of the above, we have completed the proof. \qed

\end{document}